\documentclass[12pt]{article}
\usepackage{amsmath, amsthm, amssymb,amscd}
\title{Existence of a multiplicative basis for a finitely spaced
module over an~aggregate
\footnotetext{This is the authors' version of a work that was published in Ukrainian Math. J. 46 (no. 5) (1994) 567--579.}}
\author{Andrej V. Roiter%
\\ Institute
of Mathematics\\ Tereshchenkivska 3,
Kiev, Ukraine
\and
Vladimir V. Sergeichuk%
\\ Institute
of Mathematics\\ Tereshchenkivska 3,
Kiev, Ukraine\\sergeich@ukrpack.net}
\date{}

\begin{document}
\maketitle

\noindent
It was proved in [1] that a finite-dimensional algebra, having finitely many
isoclasses of indecomposable representations, admits a multiplicative basis.
In [2] (Sections 4.10-4.12) an analogous hypothesis was formulated for
finitely spaced modules over an aggregate and an approach to its proof was
proposed. Our objective is to prove this hypothesis. Throughout this paper, $%
k$ denotes an algebraically closed field.

Let us recall some definitions from [2] (see also [3]).

By definition, an $aggregate~{\cal A}$ over $k$ is a category that satisfies
the following conditions:

a. For each $X,Y\in {\cal A}$, the set ${\cal A}(X,Y)$ is a
finite-dimensional vector space over$~k$;

b. The composition maps are bilinear;

c. ${\cal A}$ has finite direct sums;

d. Each idempotent $e \in {\cal A}(X,X)$ has the kernel.

As a consequence, each $X \in {\cal A}$ is a finite sum of indecomposables
and the algebra of endomorphisms of each indecomposable is local.

We denote by ${\cal JA}$ a spectroid of ${\cal A}$, i.e. a full subcategory
formed by chosen representatives of the isoclasses of indecomposables, and
let ${\cal R_A}$ be the radical of ${\cal A}$. We suppose that ${\cal JA}$
has finitely many objects. For each $a,b\in {\cal JA}$, the space ${\cal R}_{%
{\cal A}}(a,b)$ consists of all irreversible morphisms of ${\cal A}(a,b)$,
therefore, ${\cal A}(a,b)={\cal R}_{{\cal A}}(a,b)$ for $a\not =b$. ${\cal A}%
(a,a)=k\;$1$_a\oplus _k{\cal R}_{{\cal A}}(a,a)$.

A $module$ $M$ over an aggregate ${\cal A}$ consists of finite-dimensional
vector spaces $M(X)$, one for each object $X\in {\cal A},$ and of linear
maps $M(f):M(X)\rightarrow M(Y),$ $m\mapsto fm,\quad f\in {\cal A}(X,Y),$
which satisfy the standard axioms: $1_Xm=m,(f+g)m=fm+gm,(gf)m=g(fm),f(\alpha
m)=\alpha (fm)=(\alpha f)m,\alpha \in k.$ It gives a $k$-linear functor from
${\cal A}$ into the category mod $k$ of finite-dimensional vector spaces
over $k$. A module $M$ over ${\cal A}$ is $faithful$ if $M(f)\not =0$ for
each nonzero $f\in {\cal A}(X,Y).$

Define the $basis$ of $(M,{\cal A})$ as a set $\{m_i^a,f_l^{ba}\}$
consisting of bases $m_1^a,m_2^a,...$ of the spaces $M(a),a\in {\cal JA},$
and bases $f_1^{ba},f_2^{ba},...$ of the spaces ${\cal R_A}(a,b),a,b\in
{\cal JA}.$ The maximal rank of $M(f_l^{ba})$ is called the $rank$ of a
basis. A basis is called a $scalarly$ $multiplicative$ $basis$ if it
satisfies the following conditions:

a) Each morphism $f_l^{ba}$ is $thin$, i.e. $f_l^{ba}=g+h$ implies rank$%
~M(f_l^{ba})\leq $ rank$~M(g)$ and rank $M(f_l^{ba})\leq $rank $M(h)$ for all
$g,h\in {\cal A}(a,b);$

b) Each product $f_l^{ba}m_i^a$ has the form $\lambda m_p^b,\lambda \in k;$

c) $f_l^{ba}m_i^a=\lambda m_p^b$, $f_l^{ba}m_j^a=\mu m_p^b$, and $\lambda
,\mu \in k\setminus \{0\}$ imply $i=j.$

We say that the basis is $multiplicative$ if each nonzero product $%
f_l^{ba}m_i^a$ is a basis vector $m_p^b.$

We denote by $M^k$ the aggregate formed by all triples $(V,h,X),$ where $%
V\in $~mod$~k$, $X\in {\cal A},$ and $h\in $ Hom$_k(V,M(X)).$ A morphism
from $(V,h,X)$ to $(V^{\prime },h^{\prime },X^{\prime })$ is defined by the
pair of morphisms $\varphi \in $Hom$_k(V,V^{\prime })$ and $\xi \in {\cal A}%
(X,X^{\prime })$ such that $h^{\prime }\varphi =M(\xi )h.$ We call these
triples $spaces$ on $M$. We say that $M$ is $finitely$ $spaced$ if $M^k$ has
a finite spectroid.

The objective of the paper is to prove the following theorem:

{\bf Theorem.} {\it If }$M$ {\it is a faithful finitely spaced module over
an aggregate }${\cal A,}$ {\it then} $(M,{\cal A})$ {\it admits a
multiplicative basis of} rank$~\leq $ 2.

We wish to express gratitude to P.Gabriel, Th.Br\"ustle, T.Guidon and
U.Hassler for discussions and essential corrections.\bigskip\

{\bf 1. Construction of a scalarly multiplicative basis.}\\In Sections 1-3, $%
M$ always denotes a finitely spaced module over an aggregate$~{\cal A}$.

As shown in [2] (sections 4.7, 4.8), for each $a\in {\cal JA},$ the space $%
M(a)$ has a dimension $d(a)\leq 3$ and a sequence $m_1,m_2,...,m_{d(a)},$
where%
$$
m_i\in ({\cal R_A}(a,a))^{i-1}M(a)\setminus ({\cal R_A}(a,a))^iM(a),
$$
is a basis of $M(a)$. It will be called a $triangular$ $basis$ because the
matrix of each map $M(f),f\in {\cal A}(a,a),$ has a lower triangular form.
We assume that each basis $m_1^a,...,m_{d(a)}^a$ in a scalarly
multiplicative basis is triangular (it is always triangular up to
permutations of vectors).

A scalarly multiplicative basis is called $normed$ if it satisfies the
following condition:

d) $f_l^{ba}m_i^a=\lambda m_p^b$ and $\lambda \notin \{0,1\}$ imply that $%
f_l^{ba}m_{i^{\prime }}^a=\lambda m_{p^{\prime }}^b$ for some $i^{\prime
}<i. $

A scalarly multiplicative basis can be reduced to a normed basis by means of
multiplication of $f_l^{ba}$ by scalars.

A scalarly multiplicative basis is called $reduced$ if it satisfies
condition d) and the following condition:

e) if a morphism $\varphi =\sum_l\lambda _lf_l^{ba}$ is a product of basis
morphisms, then
$$
{\rm rank~}M(\varphi )=\sum_{\lambda_l \not =0}{\rm rank\ }M(f_l^{ba}).
$$

At the end of this section, we shall prove that every multiplicative basis
of $(M,{\cal A})$ is reduced if char $k\not =$ 2.

Let $m_1^a,...,m_{d(a)}^a$ be a fixed triangular basis of $M(a)$ for each $%
a\in {\cal JA}.$ For $m_j^a$ and $m_i^b,$ we define a linear map $%
e_{ij}^{ba}:M(a)\rightarrow M(b)$ such that $e_{ij}^{ba}m_j^a=m_i^b$ and $%
e_{ij}^{ba}m_{j^{\prime }}^a=0$ for all $j^{\prime }\not =j.$

Let $f\in {\cal R_A}(a,b),a,b\in {\cal JA.}$ We say that $f$ is a $short$ $%
morphism$ if $f\notin {\cal R_A}(c,b){\cal R_A}(a,c)$ for all $c\in {\cal %
JA,\ }f$ is a $prime$ $morphism$ if $M(f)=e_{ij}^{ba},$ and $f$ is a $double$
$morphism$ if%
$$
M(f)=e_{ij}^{ba}+\lambda e_{i^{\prime }j^{\prime }}^{ba},e_{ij}^{ba}\notin
M(a,b),i<i^{\prime },j<j^{\prime },0\not =\lambda \in k.
$$

The coefficient $\lambda $ is called the $parameter$ of a double morphism.

{\bf Proposition 1.} {\it A set} $\{m_i^a,f_l^{ba}\}$ {\it is a normed
(reduced, respectively) scalarly multiplicative basis if and only if the
following conditions are satisfied:}

1) $m_1^a,m_2^a,...${\it \ is a triangular basis of} $M(a),a\in {\cal JA,}$

2) $f_1^{ba},f_2^{ba},...${\it \ is the set of all prime and double
morphisms of} ${\cal A}(a,b),$ $a,b\in {\cal JA,}$ {\it except a single
double morphism (a single short double morphism, respectively) if the number
of double morphisms is equal to 3. Moreover, the number of double morphisms
of }${\cal A}(a,b)${\it \ is equal 0, 1 or 3, and, in the last case, there
exists a short double morphism.}

The statement of Proposition 1 about a normed scalarly multiplicative basis
follows from Lemmas 1 and 5. The complete proof of Proposition 1 will be
given in Section 3.

{\bf Lemma 1}. {\it If }$d(a)=$ 2, {\it then }$%
M(a,a)=k1_{M(a)}+ke_{21}^{aa}. $ {\it If }$d(a)=3,$ {\it then }$%
M(a,a)=k1_{M(a)}+ke_{21}^{aa}+ke_{32}^{aa}$ {\it or}
\begin{equation}
\label{d1}M(a,a)=k1_{M(a)}+k(e_{21}^{aa}+\lambda
_{aa}e_{32}^{aa})+ke_{31}^{aa}
\end{equation}
{\it and }0$\not =\lambda _{aa}\in k.$

The proof of Lemma 1 is obvious.

For every linear map $\varphi :M(a)\rightarrow M(b),$ we denote by $\varphi
_{ij}\in ke_{ij}^{ba}$ linear maps such that $\varphi =\sum \varphi _{ij}.$
We introduce an order relation on $\{1,2,...,d(b)\}\times \{1,2,...,d(a)\}$
by $(i,j)\geq (l,r)$ if $i\leq l$ and $j\geq r.$ A pair $(l,r)$ is called a $%
step$ of $\varphi \in M(a,b)$ if $\varphi _{lr}\not =0$ and $\varphi _{ij}=0$
for all $(i,j)>(l,r).$ A pair $(l,r)$ is called a {\it step} of $M(a,b)$ if $%
\psi _{lr}\not =0$ for some $\psi \in M(a,b)$ and $\varphi _{ij}=0$ for all $%
\varphi \in M(a,b)$ and all $(i,j)>(l,r)$ $(l\geq r$ because each basis $%
m_1^a,m_2^a,...$ is a triangular).

{\bf Lemma 2.} {\it If }$a,b\in {\cal JA,}a\not =b,d(a)=d(b)=3,$ {\it and } $%
M(a,b)$ {\it has two steps (1,2) and (2,3), then }$M(b,a)=ke_{31}^{ab}.$

{\bf Proof.} Let $\psi \in M(b,a).$ There is $\varphi \in M(a,b)$ having the
steps (1,2) and (2,3). By Lemma 1, there exist $\epsilon \in M(a,a)$ and $%
\delta \in M(b,b)$ such that $\varphi ^{\prime }=\varphi \epsilon +\delta
\varphi $ has the steps (1,1), (2,2) and (3,3). The inclusion ${\cal A}(b,a)%
{\cal A}(a,b)\subset {\cal R_A}(a,a)$ implies $M(b,a)M(a,b)\subset M({\cal %
R_A}(a,a))=ke_{21}^{aa}+ke_{31}^{aa}+ke_{32}^{aa}.$ Since $\psi \varphi
^{\prime }\in M({\cal R_A}(a,a)),$ all steps of $\psi $ are not higher that
(2,1) and (3,2). Since $\psi \varphi \in M({\cal R_A}(a,a)),$ we have $\psi
\in ke_{31}^{ab}.$

Therefore, $M(b,a)\subset ke_{31}^{ab}.$ Assume that $M(b,a)=0.$ Let us
examine the space ${\cal H}_\lambda =(k^6,h_\lambda ,a^2\oplus b^2)\in M^k,$
where $k^6=k\oplus k\oplus k\oplus k\oplus k\oplus k,$ $a^2=a\oplus a,$ $%
b^2=b\oplus b,$ $\lambda \in k,$ and h$_\lambda $ is the linear mapping of $%
k^6$ into $M(a^2\oplus b^2)=(km_1^a)^2\oplus (km_2^a)^2\oplus
(km_3^a)^2\oplus (km_1^b)^2\oplus (km_2^b)^2\oplus (km_3^b)^2$ with the
matrix%
$$
\left(
\begin{array}{c}
1 \\
0
\end{array}
\right) \oplus \left(
\begin{array}{c}
0 \\
1
\end{array}
\right) \oplus \left(
\begin{array}{cc}
1 & 0 \\
0 & 1
\end{array}
\right| \left.
\begin{array}{cc}
1 & 0 \\
0 & 1
\end{array}
\right) ^T\oplus \left(
\begin{array}{c}
1 \\
1
\end{array}
\right) \oplus \left(
\begin{array}{c}
1 \\
\lambda
\end{array}
\right) .
$$

We show that ${\cal H}_\lambda \not \simeq {\cal H}_\mu $ if $\lambda \not
=\mu $. Let $(\varphi ,\xi )$ be an isomorphism ${\cal H}_\lambda
\rightarrow {\cal H}_\mu .$ The linear mapping $M(\xi )$ has the block
matrix $(K_{ij}),$ $i,j\leq 6,$ where $K_{ij}$ are $2\times 2$-matrices. By $%
M(b,a)=0$ and Lemma 1, we have $K_{ij}=0$ if $i<j.$ Evidently, $%
K_{11}=K_{22}=K_{33},$ $K_{44}=K_{55}=K_{66}$ and $K_{43}=0.$

Since $h_\mu \varphi =M(\xi )h_\lambda ,$ the matrix of the nondegenerate
mapping $\varphi $ also has the block form $(\Phi _{ij}),i,j\leq 5$, where
the blocks $\Phi _{11},$ $\Phi _{22},$ $\Phi _{44}$ and $\Phi _{55}$ are 1$%
\times 1$-matrices, the block $\Phi _{33}$ is a 2$\times 2$-matrix, and $%
\Phi _{ij}=0$ if $i<j.$ Moreover,
$$
\left(
\begin{array}{c}
1 \\
0
\end{array}
\right) \Phi _{11}=K_{11}\left(
\begin{array}{c}
1 \\
0
\end{array}
\right) ,\qquad \left(
\begin{array}{c}
0 \\
1
\end{array}
\right) \Phi _{22}=K_{22}\left(
\begin{array}{c}
0 \\
1
\end{array}
\right) ,
$$
$$
\left(
\begin{array}{cc}
1 & 0 \\
0 & 1
\end{array}
\right| \left.
\begin{array}{cc}
1 & 0 \\
0 & 1
\end{array}
\right) ^T\Phi _{33}=(K_{33}\oplus K_{44})\left(
\begin{array}{cc}
1 & 0 \\
0 & 1
\end{array}
\right| \left.
\begin{array}{cc}
1 & 0 \\
0 & 1
\end{array}
\right) ^T,
$$
$$
\left(
\begin{array}{c}
1 \\
1
\end{array}
\right) \Phi _{44}=K_{55}\left(
\begin{array}{c}
1 \\
1
\end{array}
\right) ,\qquad \left(
\begin{array}{c}
1 \\
\mu
\end{array}
\right) \Phi _{55}=K_{66}\left(
\begin{array}{c}
1 \\
\lambda
\end{array}
\right) .
$$
By the third equality, we obtain $K_{33}=K_{44},$ by the first and second
equalities, we get%
$$
K_{11}=K_{22}=...=K_{66}=\left(
\begin{array}{cc}
\alpha & 0 \\
0 & \beta
\end{array}
\right) ,
$$
and, by the forth and fifth equalities, $\alpha =\beta $ and $\lambda =\mu .$
We have infinitely many nonisomorphic indecomposable spaces ${\cal H}%
_\lambda ,\lambda \in k,$ on $M.$ This proves Lemma 2.

Let $(l_1,r_1),...,(l_t,r_t)$ be all steps of $M(a,b).$ Set $%
S(a,b)=\sum_{(i,j)}ke_{ij}^{ba}$ (resp. $\overline{S}(a,b)=%
\sum_{(i,j)}ke_{ij}^{ba}$), where the sum is taken over all $(i,j)$ such
that there exists a step $(l_p,r_p)>(i,j)$ ($(l_p,r_p)\geq (i,j),$
respectively).

{\bf Lemma 3}. {\it Let} $a\not =b$\ {\it and} $M(a,b)$\ {\it have the steps}
$(1,1)$, $(2,2)${\it \ and} $(3,3)$. {\it Then there is no} $\psi \in M(a,b)$%
{\it \ such that }$M(a,b)=k\psi +S(a,b).$

{\bf Proof.} Assume that there exists $\psi \in M(a,b)$ such that $%
M(a,b)=k\psi +S(a,b).$ By the form of $M(a,b)$ and ${\cal A}(b,a){\cal A}%
(a,b)\subset {\cal R_A}(a,a),$ we have $M(b,a)\subset
ke_{21}^{ab}+ke_{31}^{ab}+ke_{32}^{ab}.$

Let us examine the space ${\cal H}_\lambda =(k^3,h_\lambda ,a\oplus b),$
where $\lambda \in k$ and $h_\lambda $ is the linear map from $k^3$ into%
$$
M(a\oplus b)=km_1^a\oplus km_2^a\oplus km_3^a\oplus km_1^b\oplus
km_2^b\oplus km_3^b
$$
with the matrix%
$$
\left(
\begin{array}{ccc}
1 & 0 & 0 \\
0 & 1 & 0 \\
0 & 0 & 1
\end{array}
\right| \left.
\begin{array}{ccc}
0 & 0 & 0 \\
0 & 1 & 0 \\
0 & 0 & \lambda
\end{array}
\right) ^T.
$$

Let $(\varphi ,\xi )$ be an isomorphism ${\cal H}_\lambda \rightarrow {\cal H%
}_\mu .$ It follows from the conditions imposed on $M(a,a),M(a,b),M(b,a)$
and $M(b,b)$ that the matrix of $M(\xi )$ has form

$$
\left(
\begin{array}{ccc|ccc}
\alpha _1 & 0 & 0 &0 & 0 & 0 \\
\alpha _2 & \alpha _1 & 0 &\gamma _1 & 0 & 0  \\
\alpha _4 & \alpha _3 & \alpha _1 & \gamma _3 & \gamma _2 & 0  \\ \hline
\delta _1 & 0 & 0 &
 \beta _1 & 0 & 0
  \\
\delta _4 & \delta _2 & 0 &
\beta _2 & \beta _1 & 0
\\
\delta _6 & \delta _5 & \delta _3&\beta _4 & \beta _3 & \beta _1
\end{array}
\right) .
$$
Moreover, $\delta _1=\delta \epsilon _1$, $\delta _2=\delta \epsilon _2$,
and $\delta _3=\delta \epsilon _3$, where $\delta \in k$ and $\epsilon _1$,$%
\epsilon _2$, and $\epsilon _3$ are the diagonal elements of the lower
triangular matrix of $\psi .$ By $h_\mu \varphi =M(\xi )h_\lambda ,$ we find
successively that $\delta =0$, the mapping $\varphi $ has the lower
triangular matrix with the diagonal $(\alpha _1,\alpha _1,\alpha _1),$ $%
\alpha _1=\beta _1,$ and $\lambda =\mu .$

Hence ${\cal H}_\lambda \not \simeq {\cal H}_\mu $ for $\lambda \not =\mu $
and $M$ is infinitely spaced. We arrive at a contradiction that proves
Lemma~3.

{\bf Lemma 4}. $S(a,b)\subset M(a,b).$

{\bf Proof.} We must show that if $(l,r)$ is a step of $M(a,b)$, then%
$$
S_{lr}(a,b)=\sum_{(i,j)<(l,r)}ke_{ij}^{ba}\subset M(a,b).
$$
By Lemma 3, there exists a $\psi \in M(a,b)$ having the step $(l,r)$ but not
more then two steps. If $\psi $ and $M(a,b)$ have the steps $(1,2)$ and $%
(2,3)$, then, by Lemma 2, $e_{31}^{ab}\psi \in M(a,a)$ has the unique step $%
(3,2)$. Hence,%
$$
M(a,a)=k1_{M(a)}\oplus ke_{21}^{aa}\oplus ke_{31}^{aa}\oplus ke_{32}^{aa}.
$$
In all other cases, by Lemma 1, $S_{lr}(a,b)$ is contained in the space
generated by all $\delta \psi \epsilon ,$ where $\epsilon \in M(a,a)$ and $%
\delta \in M(b,b)$. This proves Lemma 4.

By Lemma 4, we have the following lemma.

{\bf Lemma 5}. {\it Let} $a,b\in {\cal JA,}a\not =b,$ {\it and} $M(a,b)\not
= \overline{S}(a,b).$ {\it Then} {\it only} {\it three} {\it cases can occur}
$(\lambda _{ab}\not =0\not =\mu _{ab}):$

{\it a)} $M(a,b)$ {\it has two steps} $(l_1,r_1)$ {\it and} $(l_2,r_2),$ $%
l_1<l_2,$ {\it and is equal to}%
$$
k(e_{l_1r_1}^{ba}+\lambda _{ab}e_{l_2r_2}^{ba})\oplus S(a,b);
$$

{\it b)} $M(a,b)$ {\it has the steps} $(1,1),(2,2)$ {\it and} $(3,3)$ {\it %
and is equal to}%
$$
k(e_{11}^{ba}+\lambda _{ab}e_{22}^{ba})\oplus ke_{33}^{ba}\oplus S(a,b),
$$
{\it or}%
$$
k(e_{11}^{ba}+\lambda _{ab}e_{33}^{ba})\oplus ke_{22}^{ba}\oplus S(a,b),
$$
{\it or}%
$$
k(e_{22}^{ba}+\lambda _{ab}e_{33}^{ba})\oplus ke_{11}^{ba}\oplus S(a,b);
$$

{\it c)} $M(a,b)$ {\it has the steps} $(1,1),(2,2)$ {\it and} $(3,3)$ {\it %
and is equal to}%
$$
k(e_{11}^{ba}+\lambda _{ab}e_{22}^{ba})\oplus k(e_{11}^{ba}+\mu
_{ab}e_{33}^{ba})\oplus S(a,b).
$$

{\bf Remarks.} 1) In a normed scalarly multiplicative basis, each long
double morphism $\varphi \in {\cal A}(a,b)$ is the product of double basis
morphisms. Indeed, let $\varphi =\tau \psi $, where $\psi \in {\cal R_A}%
(a,c) $ and $\tau \in {\cal R_A}(c,b)$. Then $\psi $ is the unique double
morphism of ${\cal A}(a,c)$ (otherwise, $\varphi $ is the sum of prime
morphisms). Therefore, $\psi $ is a basis morphism. Similarly, $\tau $ is
also a basis morphism.

2) A normed scalarly multiplicative basis is reduced if and only if all long
double morphisms are basis morphisms. Indeed, let a long double morphism $%
\varphi \in {\cal A}(a,b)$ be not a basis morphism. Then ${\cal A}(a,b)$ has
two double morphisms and $\varphi $ is their linear combination. But this
contradicts the definition of a reduced basis.

3) Lemma 1 and Lemma 5 imply the statement of Proposition 1 about a normed
scalarly multiplicative basis. By Remark 2, to complete the proof of
Proposition 1 we must prove that each ${\cal A}(a,b)$ $(a,b\in {\cal JA})$
does not contain three long double morphisms.

4) If char $k\not =2,$ then every multiplicative basis is reduced. Indeed,
otherwise, there is, by Remark 2, a long double morphism $\varphi \in {\cal
A}(a,b)$, which is not a basis morphism. By Lemma 5, $\varphi =\psi -\tau
$, where $\psi $ and $\tau $ are basis long double morphisms of ${\cal
A}(a,b)$
$\varphi $ is a product of basis morphisms; hence $M(\varphi
)=e_{ii}^{ba}+e_{jj}^{ba}$ and char $k=2$.\bigskip\

{\bf 2. The graph of a scalarly multiplicative basis.}

In this section, we study some properties of a scalarly multiplicative basis
and give the proof of Proposition~1.

Following [2] (Section 4.9), we define a poset $~{\cal P},$ whose elements
are the spaces $a_i=~({\cal R_A}(a,a))^{i-1}M(a)$ $(a\in {\cal JA},$ $1\leq
i\leq d(a))$ and where $a_i\leq b_j$ if and only if ${\cal A}(b,b)fa_i=b_j$
for some $f\in {\cal A}(a,b).$ The elements $a_i\in {\cal P}$ are in a
one-to-one correspondence with the basis vectors $m_i^a$ of every scalarly
multiplicative basis $\{m_i^a,f_l^{ba}\}$, moreover, $a_i<b_j$ if and only
if $f_l^{ba}m_i^a=\lambda m_j^a$ for some $f_l^{ba}$ and $0\not =\lambda \in
k.$ We decompose the poset ${\cal P\ }$into disjoint totally ordered subsets
$\{a_1,...,a_{d(a)}\},$ $(a_1<a_2<...<a_{d(a)},d(a)\leq 3);$ each of them is
called a {\it double} if $d(a)=2$ and a {\it triple} if $d(a)=3.$

The following three lemmas were given in [2] without proofs.

{\bf Lemma 6} (see [2] (Lemma 4.12.1)). {\it The union} $\cup
\{a_1,a_2,a_3\} $ {\it of} {\it all} {\it triples} {\it is} {\it totally}
{\it ordered}.

{\bf Proof}. The elements of a triple are totally ordered.

Let $\{a_1,a_2,a_3\}$ and $\{b_1,b_2,b_3\}$ be triples and let some $a_i$ be
not comparable with some $b_j.$ We shall construct indecomposable spaces $%
{\cal H}_\lambda =(k^6,h_\lambda ,a^2\oplus b^2)$ on $M$, $\lambda \in k,$
such that ${\cal H}_\lambda \not \simeq {\cal H}_\mu $ for $\lambda \not
=\mu .$

For $i=3$ and $j=1$, the spaces ${\cal H}_\lambda $ were constructed in the
proof of Lemma 2. For arbitrary $i$ and $j$, ${\cal H}_\lambda $ is
constructed analogously with the block%
$$
\left(
\begin{array}{cc}
1 & 0 \\
0 & 1
\end{array}
\right| \left.
\begin{array}{cc}
1 & 0 \\
0 & 1
\end{array}
\right) ^T
$$
of $h_\lambda :k^6\rightarrow M(a^2\oplus b^2)$ located in the rows of%
$$
km_i^a\oplus km_i^a\oplus km_j^b\oplus km_j^b\subset M(a^2\oplus b^2).
$$
Let $(\varphi ,\xi ):{\cal H}_\lambda \tilde{\rightarrow }{\cal H}_\mu $ and
let $(M_{ij})$ be the block matrix of $M(\xi )$. Then $(M_{ij})$ is not
upper block-triangular, but we can reduce $(M_{ij})$ to the upper
block-triangular form by means of simultaneous transpositions of vertical
and horizontal stripes, since the set $\{a_1,a_2,a_3,b_1,b_2,b_3\}$ is
partially ordered. Hence, $M$ is infinitely spaced. We arrive at a
contradiction that proves Lemma 6.

{\bf Lemma 7} (see [2] (Lemma 4.9)). {\it There are no elements} $a_i$, $%
a_{i^{\prime }},$ $b_j,$ {\it and} $b_{j^{\prime }}$ {\it such} {\it that }$%
a_i\not =a_{i^{\prime }},$ $b_j\not =b_{j^{\prime }}$, $a_i$ {\it is not}
{\it comparable} {\it to }$b_{j^{\prime }}$, {\it and} $b_j$ {\it is} {\it %
not} {\it comparable} {\it to} $a_{i^{\prime }}$. {\it There} {\it are} {\it %
no} {\it elements }$a_i$, $a_{i^{\prime }},$ $b_j,$ $b_{j^{\prime }}$ $c_l$
{\it and} $c_{l^{\prime }}$ {\it such} {\it that }$a_i\not =a_{i^{\prime }},$
$b_j\not =b_{j^{\prime }}$, $c_l\not =c_{l^{\prime }},$ $a_i$ {\it is not}
{\it comparable} {\it to }$b_{j^{\prime }}$, $b_j$ {\it is} {\it not} {\it %
comparable} {\it to} $c_{l^{\prime }},$ {\it and} $c_l$ {\it is} {\it not}
{\it comparable} {\it to} $a_{i^{\prime }}$.

{\bf Proof.} In the first case, we set ${\cal H}_\lambda =(ke_1\oplus
ke_2,h_\lambda ,a\oplus b)\in M^k,$ where $h_\lambda e_1=m_i^a+m_{j^{\prime
}}^b$ and $h_\lambda e_2=m_j^b+\lambda m_{i^{\prime }}^a.$ In the second
case, we set ${\cal H}_\lambda =(ke_1\oplus ke_2\oplus ke_3,h_\lambda
,a\oplus b\oplus c),$ where $h_\lambda e_1=m_i^a+m_{j^{\prime }}^b,$ $%
h_\lambda e_2=m_j^b+m_{l^{\prime }}^c,$ and $h_\lambda e_3=m_l^c+\lambda
m_{i^{\prime }}^a.$ Obviously, ${\cal H}_\lambda \not \simeq {\cal H}_\mu $
for $\lambda \not =\mu .$

{\bf Lemma 8} (see [2] (Lemma 4.12.2)). {\it Each} {\it triple} {\it contains%
} {\it at} {\it least} {\it two} {\it elements} {\it comparable} {\it with}
{\it all} {\it elements of all doubles}.

{\bf Proof.} Assume that Lemma 8 is not true for a triple $\{a_1,a_2,a_3\}$
and doubles $\{b_1,b_2\}$ and $\{c_1,c_2\}.$

{\bf Case 1}. Assume that $b\not =c$. For definiteness, we suppose that $a_2$
is not comparable to $b_1$ and $a_3$ is not comparable to $c_1$.

For each representation ${\cal H}$\\%
\begin{picture}(0,0)
\put(150,0){$k^{r_1}$}
\put(137,-15){$A_1$}
\put(150,-15){$\downarrow$}

\put(122,-23){$A_2$}
\put(90,-23){$B_1$}
\put(52,-23){$B_2$}
\put(170,-23){$A_3$}
\put(205,-23){$C_1$}
\put(245,-23){$C_2$}

\put(150,-30){$k^{t_1}$}
\put(120,-30){$\longrightarrow$}
\put(105,-30){$k^{r_2}$}

\put(85,-30){$\longleftarrow$}
\put(70,-30){$k^{t_2}$}
\put(50,-30){$\longrightarrow$}
\put(35,-30){$k^{r_4}$}
\put(165,-30){$\longleftarrow$}
\put(188,-30){$k^{r_3}$}
\put(203,-30){$\longrightarrow$}
\put(225,-30){$k^{t_3}$}
\put(240,-30){$\longleftarrow$}
\put(260,-30){$k^{r_5}$}

\end{picture}

\bigskip\bigskip\ \\of the quiver $\tilde E_7$ (see [2], (Section 6.3)), we
construct the space%
$$
\overline{{\cal H}}=(k^{r_1+...+r_5},h,a^{t_1}\oplus b^{t_2}\oplus
c^{t_3})\in M^k,
$$
where%
$$
h=A_1\oplus \left(
\begin{array}{c}
A_2 \\
B_1
\end{array}
\right) \oplus \left(
\begin{array}{c}
A_3 \\
C_1
\end{array}
\right) \oplus B_2\oplus C_2
$$
is a linear mapping of $k^{r_1+...+r_5}$ into
$$
\begin{array}{c}
M(a^{t_1}\oplus b^{t_2}\oplus c^{t_3})=(km_1^a)^{t_1}\oplus \left[
(km_2^a)^{t_1}\oplus (km_1^b)^{t_2}\right] \oplus \\
\left[ (km_3^a)^{t_1}\oplus (km_1^c)^{t_3}\right] \oplus
(km_2^b)^{t_2}\oplus (km_2^c)^{t_3}.
\end{array}
$$

The functor ${\cal H\mapsto }\overline{{\cal H}}$ on the representations $%
{\cal H}$ with injective $A_1,$ $A_2,$ $A_3,$ $B_2,$ and $C_2$ preserves
indecomposability and {\it heteromorphism} (i.e. ${\cal H\simeq H^{\prime }}$
if $\overline{{\cal H}}\simeq {\cal \overline{H}^{\prime }}$). Indeed, let $%
(\varphi ,\xi ):\overline{{\cal H}}\tilde \rightarrow \overline{{\cal H}}%
{\cal ^{\prime }.}$ The nondegenerate linear maps $\varphi $ and $M(\xi )$
have the block forms $(\Phi _{ij}),$ $i,j\leq 5,$ and $(K_{ij}),$ $i,j\leq 7.
$ The equality $h^{\prime }\varphi =M(\xi )h$ implies $A_1^{\prime }\Phi
_{11}=K_{11}A_1,$
$$
\begin{array}{c}
\left(
\begin{array}{c}
A_2^{\prime } \\
B_1^{\prime }
\end{array}
\right) \Phi _{22}=\left(
\begin{array}{cc}
K_{22} & K_{23} \\
K_{32} & K_{33}
\end{array}
\right) \left(
\begin{array}{c}
A_2 \\
B_1
\end{array}
\right) ,\qquad B_2^{\prime }\Phi _{44}=K_{66}B_2, \\
\\
\left(
\begin{array}{c}
A_3^{\prime } \\
C_1^{\prime }
\end{array}
\right) \Phi _{33}=\left(
\begin{array}{cc}
K_{44} & K_{45} \\
K_{54} & K_{55}
\end{array}
\right) \left(
\begin{array}{c}
A_3 \\
C_1
\end{array}
\right) ,\qquad C_2^{\prime }\Phi _{55}=K_{77}C_2.
\end{array}
$$

Since $\{a_1,a_2,a_3\}$ is a triple and $\{b_1,b_2\}$ and $\{c_1,c_2\}$ are
doubles, we have $K_{11}=K_{22}=K_{44},$ $K_{33}=K_{66},$ and $K_{55}=K_{77}$%
. Since $a_2$ is not comparable to $b_1$ and $a_3$ is not comparable to $c_1,
$ we have $K_{23}=0,$ K$_{32}=0,$ $K_{45}=0,$ and $K_{54}=0.$ Hence, the
diagonal blocks of $(\Phi _{ij})$ and $(K_{ij})$ determine a morphism ${\cal %
H\rightarrow H^{\prime }.}$

We shall show that this morphism is an isomorphism, i.e. the diagonal blocks
$\Phi _{ii}$ and $K_{ii}$ are invertible. By strengthening the partial order
relation in $\{a_1,$ $a_2,$ $a_3,$ $b_1,$ $b_2,$ $c_1,$ $c_2\}$, we obtain a
total order relation $\ll $ such that $a_2\ll b_1$ and $a_3\ll c_1$ (these
pairs are not comparable with respect to $<)$.

We transpose the horizontal stripes of the matrices of $h$ and $h^{\prime }$
according to the new order. Then we transpose the vertical stripes to get
lower trapezoidal matrices. Correspondingly, we transpose the blocks of $%
(\Phi _{ij})$ and $(K_{ij}).$ Then the new matrix $(K_{ij})$ has a lower
triangular form. The upper nonzero blocks of vertical stripes are the
injective maps $A_1,$ $A_2,$ $A_3,$ $B_2,$ and $C_2$ (since $a_2\ll b_1$ and
$a_3\ll c_1).$ It follows from $h^{\prime }\varphi =M(\xi )h$ that $(\Phi
_{ij})$ also has a lower triangular form. Hence, the diagonal blocks $\Phi
_{ii}$ and $K_{ii}$ are invertible and ${\cal H\simeq H^{\prime }.}$

But the quiver $\tilde E_7$ admits an infinite set of nonisomorphic
indecomposable representations of the form ${\cal H}$ with injective $A_1,$ $%
A_2,$ $A_3,$ $B_2,$ and $C_2$ (and surjective $B_1$ and $C_1,$ which will be
used in the case 2). These representations are determined by the matrices%
$$
\begin{array}{c}
(A_1\mid A_2\mid A_3)=\left(
\begin{array}{cc}
1 & \alpha  \\
1 & 1 \\
1 & 0 \\
0 & 1
\end{array}
\right| \left.
\begin{array}{ccc}
1 & 0 & 0 \\
0 & 1 & 0 \\
0 & 0 & 1 \\
0 & 0 & 0
\end{array}
\right| \left.
\begin{array}{ccc}
0 & 0 & 0 \\
0 & 0 & 1 \\
0 & 1 & 0 \\
1 & 0 & 0
\end{array}
\right) , \\
\\
(B_1\mid B_2)=(C_1\mid C_2)=\left(
\begin{array}{ccc}
0 & 1 & 0 \\
0 & 0 & 1
\end{array}
\right| \left.
\begin{array}{c}
1 \\
0
\end{array}
\right) ,
\end{array}
$$
and they are nonisomorphic for different $\alpha \in k.$ This contradicts
the assumption that $M$ is finitely spaced.

{\bf Case 2}. Assume that $b=c$. By Lemma 7, if $a_i$ is not comparable to $%
b_1$, and $a_j$ is not comparable to $b_2,$ then $i=j$. Let $a_2$ and $a_3$
be not comparable to $b_1.$ Then $a_1<b_1$ and $a_3<b_2$.

As in the case 1, for each representation ${\cal H}$ of the quiver $\tilde{E}%
_7$ with injective $A_1,$ $A_2,$ $A_3,$ $B_2,$ and $C_2$ and surjective $B_1$
and $C_1$, we construct the space $\widehat{{\cal H}}%
=(k^{r_1+...+r_5},h,a^{t_1}\oplus b^{t_2+t_3})\in M^k$, where%
$$
h=A_1\oplus \left(
\begin{array}{c}
A_2 \\
B_1
\end{array}
\right) \oplus \left(
\begin{array}{c}
A_3 \\
C_1
\end{array}
\right) \oplus B_2\oplus C_2
$$
is a linear mapping of $k^{r_1+...+r_5}$ into
$$
\begin{array}{c}
M(a^{t_1}\oplus b^{t_2+t_3})=(km_1^a)^{t_1}\oplus \left[
(km_2^a)^{t_1}\oplus (km_1^b)^{t_2}\right] \oplus \\
\left[ (km_3^a)^{t_1}\oplus (km_1^b)^{t_3}\right] \oplus
(km_2^b)^{t_2}\oplus (km_2^b)^{t_3}.
\end{array}
$$

Let $(\varphi ,\xi ):\widehat{\cal H}\tilde \rightarrow \widehat{\cal H}%
{\cal ^{\prime }.}$ It follows from the order relation for $%
\{a_1,a_2,a_3,b_1,b_2\}$ that all blocks over the diagonal of the block
matrix $K=(K_{ij})_{i,j=1,2,...,7}$ of the mapping $M(\xi )$ are zero except
the blocks $K_{35}=K_{67}.$ Let us prove that they are zero, too.

Indeed, by comparing the blocks with index $(2,3)$ in the equality $%
h^{\prime }\varphi =M(\xi )h,$ we obtain $A_2^{\prime }\Phi _{23}=0$ and $%
\Phi _{23}=0$ since $A_2^{\prime }$ is injective. By comparing the blocks
with index $(3,3),$ we obtain B$_1^{\prime }\Phi _{23}=K_{35}C_1$ and K$%
_{35}=0$ since C$_1$ is surjective.

Hence $K$ is the lower block-triangular matrix. Therefore $\Phi $ also is a
lower block-triangular matrix, the diagonal blocks $K_{ii}$ and $\Phi _{ii}$
of which are invertible, ${\cal H\simeq H^{\prime }}$. This proves our lemma.

Now fix a normed scalarly multiplicative basis $\{m_i^a,f_l^{ba}\}$ and
define the oriented graph $\Gamma $, the set of vertices $\Gamma _0$ of
which is the poset ${\cal P}$ and there is an arrow $a_p\rightarrow b_q$ $%
(a_p,b_q\in \Gamma _0)$ if and only if $M(f_l^{ba})=\lambda e_{qp}^{ba}+\mu
e_{q^{\prime }p^{\prime }}^{ba}$ for some short double morphism $f_l^{ba}$
(then there is an arrow $a_{p^{\prime }}\rightarrow b_{q^{\prime }}$ and we
shall say that the arrows $a_p\rightarrow b_q$ and $a_{p^{\prime
}}\rightarrow b_{q^{\prime }}$ are {\it connected}). An arrow $%
a_p\rightarrow b_q$ will be called a {\it weak arrow} if ${\cal A}(a,b)$
contains three double morphisms. Each weak arrow is connected with two
arrows. The others will be called {\it strong arrows}, each of them is
connected exactly with one arrow.

{\bf Lemma 9.} {\it Let} $a_i<b_j<c_r$ {\it and} $a_i\rightarrow c_r$ {\it be%
} {\it an} {\it arrow.} {\it Then} $a\not =b\not =c\not =a,$ $i=r,$ {\it the}
{\it spaces} ${\cal A}(a,b)$, ${\cal A}(b,c)$ {\it and} ${\cal A}(a,c)$ {\it %
contain} {\it exactly} $1,1$ {\it and} $3$ {\it double} {\it morphisms} {\it %
respectively,} {\it and} {\it there} {\it exists} {\it a pair} {\it of
oriented} {\it paths} $(a_i\rightarrow ...\rightarrow b_j\rightarrow
...\rightarrow c_i,$ $a_{i^{\prime }}\rightarrow ...\rightarrow b_{j^{\prime
}}\rightarrow ...\rightarrow c_{i^{\prime }})$ {\it consisting} {\it of}
{\it connected} {\it strong} {\it arrows,} {\it and} {\it a pair} {\it of}
{\it connected} {\it weak} {\it arrows} $(a_i\rightarrow c_i,$ $a_{i^{\prime
\prime }}\rightarrow c_{i^{\prime \prime }},$ $i\not =i^{\prime \prime }).$
{\it In} {\it the} {\it case} {\it of} {\it a reduced} {\it scalarly} {\it %
multiplicative} {\it basis,} {\it there} {\it is} {\it no} {\it other} {\it %
arrow} {\it from} $\{a_l\}$ {\it to} $\{c_l\}.$

{\bf Proof}. Since $a_i<b_j<c_r$, there are morphisms $g\in {\cal A}(a,b)$
and $h\in {\cal A}(b,c)$ such that $M(g)=\alpha e_{ji}^{ba}+\beta
e_{j^{\prime }i^{\prime }}^{ba}$ and $M(h)=\gamma e_{rj}^{cb}+\delta
e_{r^{\prime \prime }j^{\prime \prime }}^{cb}$ $(\alpha ,\beta ,\gamma
,\delta \in k$ and $\alpha \not =0\not =\gamma )$. If $hg$ is a prime
morphism, then $M(hg)=\alpha \gamma e_{ri}^{ca}$ contradicts the existence
of the arrow $a_i\rightarrow c_r$. Hence $hg$ is a double morphism, $\beta
\not =0\not =\delta ,$ $j^{\prime }=j^{\prime \prime }$ and $g$ and $h$ are
the unique double morphisms of ${\cal A}(a,b)$ and ${\cal A}(b,c)$
respectively. The space ${\cal A}(a,c)$ contains the double morphism $hg$
and the short double morphism corresponding to the arrow $a_i\rightarrow
c_r, $ hence $M(a,c)$ has the form from item c) of Lemma 5.

If the basis is reduced then by Remark 2 of Sect.1, the double morphism $hg$
is a basis morphism and there is only one pair of connected arrows from $%
\{a_l\}$ to $\{c_l\}$. This proves our lemma.

{\bf Proof of Proposition 1.} By Remark 3 of Sect.1, we must prove that each
space ${\cal A}(a,c)$ $(a,c\in {\cal JA})$ does not contain three long
double morphisms. By contradiction let $f_1,f_2,f_3\in {\cal A}(a,c)$ be
three long double morphisms and let $f_r=h_rg_r,$ where g$_r$ is a short
double morphism and $r=1,2,3.$ The morphisms $g_1,g_2$ and $g_3$ correspond
to the pairs of connected arrows $(a_1\rightarrow x_i,$ $a_2\rightarrow
x_{i^{\prime }}),$ $(a_1\rightarrow y_j,$ $a_3\rightarrow y_{j^{\prime }}),$
and $(a_2\rightarrow z_l,$ $a_3\rightarrow z_{l^{\prime }}).$

Let $x_i<y_j.$ By putting $(a_i,b_j,c_r)=(a_1,x_i,y_j)$ in Lemma 9, we
obtain that ${\cal A}(a,y)$ contains three double morphisms. By putting $%
(a_i,b_j,c_r)=(a_1,y_j,c_1)$ in Lemma 9, we have that ${\cal A}(a,y)$
contains exactly one double morphism.

Hence $x_i$ is not comparable to $y_j$. Similarly $x_{i^{\prime }}$ is not
comparable to z$_l,$ and y$_{j^{\prime }}$ is not comparable to z$%
_{l^{\prime }}$. This contradicts Lemma 7 and proves Proposition 1.

We shall now assume that the graph $\Gamma $ is obtained from a reduced
scalarly multiplicative basis.

{\bf Lemma 10.} {\it If two arrows start from (stop at) the same vertex,}
{\it then} {\it the} {\it arrows} {\it connected} {\it with} {\it them} {\it %
start} {\it from} {\it (stop} {\it at)} {\it different vertices.}

{\bf Proof.} By contradiction, let $b_j\leftarrow a_i\rightarrow c_r$ and $%
b_{j^{\prime }}\leftarrow a_{i^{\prime }}\rightarrow c_{r^{\prime }}$ be
connected arrows. If $b_j<c_r$, then $a_i<b_j<c_r$ and, by Lemma 9, the
arrows connected with $a_i\rightarrow b_j$ and $a_i\rightarrow c_r$ must
start from different vertices, but they start from $a_{i^{\prime }}.$
Analogously, $b_{j^{\prime }}$ is not comparable to $c_{r^{\prime }}$. This
contradicts Lemma 7.

{\bf Lemma 11.} {\it There are no two arrows starting from} {\it (stopping}
{\it at)} {\it the} {\it same} {\it vertex} {\it of} {\it a} {\it double.}
{\it There} {\it are} {\it no} {\it three} {\it arrows starting from} {\it %
(stopping} {\it at)} {\it the} {\it same} {\it vertex} {\it of} {\it a} {\it %
triple.}

The proof follows from Lemma 10.

{\bf Lemma 12.} {\it There are at most two different} {\it pairs} {\it of}
{\it connected} {\it arrows starting from} {\it (stopping} {\it at)} {\it the%
} {\it same} {\it triple.}

{\bf Proof.} By contradiction, let there be three pairs of connected arrows
from a triple $\{a_1,a_2,a_3\}$ to $\{b_i\},$ $\{c_i\},$ $\{d_i\}.$ Since
there exist at most two pairs of connected arrows from a triple to a triple,
then there are no three coinciding objects among $a,$ $b,$ $c,$ $d.$ Hence
there exist five possibilities up to a permutation of $b,$ $c,$ $d$: 1) $a=b%
\not =c\not =d,$ $a\not =d$; 2) $a=b\not =c=d;$ 3) $a\not =b=d\not =c,$ $a%
\not =c;$ 4) $a,$ $b,$ $c,$ $d$ are distinct and there are two arrows $%
a_i\rightarrow b_j$ and $a_i\rightarrow c_r,$ $b_j<c_r;$ 5) $a,$ $b,$ $c,$ $%
d $ are distinct and for each pair of arrows $a_i\rightarrow x$ and $%
a_i\rightarrow y,$ the vertices $x$ and $y$ are incomparable.

By Lemmas 9-11, we have the following subgraphs of $\Gamma $ in cases $1,3$
and $4$:
\newpage
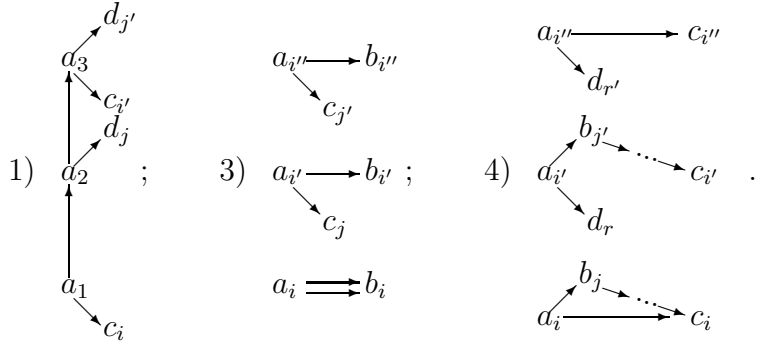
\begin{picture}(0,0)
\put(30,-63){$1)$}
\put(50,-20){$a_3$}
\put(55,-15){\vector(1,1){10}}
\put(66,-4){$d_{j^\prime}$}
\put(55,-23){\vector(1,-1){10}}
\put(66,-36){$c_{i^\prime}$}
\put(50,-63){$a_2$}
\put(55,-58){\vector(1,1){10}}
\put(66,-47){$d_{j}$}
\put(53,-57){\vector(0,1){35}}
\put(50,-106){$a_1$}
\put(55,-108){\vector(1,-1){10}}
\put(66,-122){$c_i$}
\put(53,-100){\vector(0,1){35}}
\put(80,-63){$;$}

\put(110,-63){$3)$}
\put(130,-20){$a_{i^{\prime \prime}}$}
\put(138,-23){\vector(1,-1){10}}
\put(149,-39){$c_{j^\prime}$}
\put(143,-18){\vector(1,0){20}}
\put(165,-20){$b_{i^{\prime \prime}}$}
\put(130,-63){$a_{i^{\prime}}$}
\put(138,-66){\vector(1,-1){10}}
\put(149,-82){$c_{j}$}
\put(143,-61){\vector(1,0){20}}
\put(165,-63){$b_{i^{\prime}}$}
\put(130,-106){$a_{i}$}
\put(143,-102){\vector(1,0){20}}
\put(143,-106){\vector(1,0){20}}
\put(165,-106){$b_{i}$}
\put(180,-63){$;$}

\put(210,-63){$4)$}
\put(230,-10){$a_{i^{\prime \prime}}$}
\put(238,-13){\vector(1,-1){10}}
\put(249,-29){$d_{r^\prime}$}
\put(243,-8){\vector(1,0){40}}
\put(287,-10){$c_{i^{\prime \prime}}$}
\put(230,-63){$a_{i^{\prime}}$}
\put(238,-66){\vector(1,-1){10}}
\put(249,-82){$d_r$}
\put(235,-58){\vector(1,1){10}}
\put(246,-47){$b_{j^\prime}$}
\put(288,-63){$c_{i^{\prime}}$}
\put(255,-49){\vector(3,-1){10}}
\put(267,-54){.}
\put(270,-55){.}
\put(273,-56){.}
\put(276,-56){\vector(3,-1){10}}
\put(230,-118){$a_i$}
\put(235,-113){\vector(1,1){10}}
\put(246,-102){$b_j$}
\put(288,-118){$c_i$}
\put(255,-104){\vector(3,-1){10}}
\put(267,-109){.}
\put(270,-110){.}
\put(273,-111){.}
\put(276,-111){\vector(3,-1){10}}
\put(240,-115){\vector(1,0){40}}
\put(310,-63){.}
\end{picture}

\bigskip
\bigskip
\bigskip
\bigskip
\bigskip
\bigskip
\bigskip
\bigskip
\bigskip
\bigskip
\bigskip
Consider these cases.

1) If $c_i<a_2$ or $d_j<a_3$, then by Lemma 9, ${\cal A}(a,a)$ contains
three double morphisms, which is a contradiction. If $a_2<c_i$ or $a_3<d_j,$
then by Lemma 9, there is an arrow $a_2\rightarrow c_i$ or $a_3\rightarrow
d_j,$ in contradiction with Lemma 11. Hence $a_2$ is incomparable with $c_i$
and $a_3$ is incomparable with $d_j$, which is impossible by Lemma~8.

2) This case is similar to the previous one.

3) The inequality $b_{i^{\prime }}<c_j$ is impossible, by Lemma 9, because $%
{\cal A}(a,b)$ contains three double morphisms. The inequality $b_{i^{\prime
}}>c_j$ is impossible, by Lemma 9, because there are four arrows from $%
\{a_l\}$ to $\{b_l\}$. Hence $b_{i^{\prime }}$ is incomparable with $c_j.$
Analogously $b_{i^{\prime \prime }}$ is not comparable with $c_{j^{\prime }}$
in contradiction with Lemma 7.

4) The inequalities $c_{i^{\prime }}<d_r$ and $c_{i^{\prime \prime
}}<d_{r^{\prime }}$ are impossible, by Lemma 9, because ${\cal A}(a,c)$
contains three double morphisms. If $d_r<c_{i^{\prime }}$ or $d_{r^{\prime
}}<c_{i^{\prime \prime }},$ then the double morphism $\lambda e_{i^{\prime
}i^{\prime }}^{ca}+\mu e_{i^{\prime \prime }i^{\prime \prime }}^{ca}$ ($%
\lambda \not =0\not =\mu )$ is a product of double morphisms in ${\cal A}%
(a,d)$ and ${\cal A}(d,c)$, hence ${\cal A}(a,c)$ contains two long double
morphisms in contradiction with the arrows $a_i\rightarrow c_i$ and $%
a_{i^{\prime \prime }}\rightarrow c_{i^{\prime \prime }}.$ Hence $%
c_{i^{\prime }}$ is not comparable to $d_r$ and c$_{i^{\prime \prime }}$ is
not comparable to $d_{r^{\prime }},$ in contradiction with Lemma~7.

5) This case is impossible by Lemma 7. The proof of Lemma 12 is thus
complete.\bigskip\

{\bf 3. A construction of a multiplicative basis.}

In this section we shall prove the following proposition.

{\bf Proposition 2.} {\it From every reduced scalarly multiplicative basis,
we} {\it can} {\it obtain} {\it a} {\it reduced} {\it scalarly} {\it %
multiplicative} {\it basis} {\it by means of multiplications} of {\it the}
{\it basis} {\it vectors} {\it by} {\it non-zero} {\it elements} {\it of k.}

Let $\Gamma $ be the graph of a reduced scalarly multiplicative basis $%
\{m_i^a,f_l^{ba}\}$ and let $\Gamma _1$ be the set of its arrows. An
integral function $z:\Gamma _1\rightarrow ${\cal Z} will be called a {\it %
weight} {\it function} and its value at an arrow will be called the {\it %
weight} {\it of} {\it the arrow} if:

a) $z(\alpha _1)=-z(\alpha _2)$ for each pair of connected arrows $\alpha
_1,\alpha _2$;

b) the sum of the weights of all arrows stopping at a vertex $v\in \Gamma _0$
is equal to the sum of the weights of all arrows starting from $v$ (this sum
will be called the {\it weight} {\it of} $v$ and will be denoted by $z(v)$).

{\bf Lemma 13.} {\it There} {\it exists} {\it no} {\it non-zero} {\it weight}
{\it function.}

{\bf Proof.} By contradiction let $z:\Gamma _1\rightarrow $ {\cal Z} be a
non-zero weight function. An arrow $\alpha $ will be called {\it %
nondegenerate} if $z(\alpha )\not =0$.

Let $v_1<...<v_m$ be the set of all vertices of the triples of $\Gamma .$
For each vertex $v_i,$ we denote by $v_{i^{\prime }},v_{i^{\prime \prime }}$
the two vertices such that $\{v_i,v_{i^{\prime }},v_{i^{\prime \prime }}\}$
is a triple.

By an {\it elementary path of length s} we shall mean a sequence of arrows
of the form
\begin{equation}
\label{d2}
\begin{array}{c}
\\
v_p
\end{array}
\begin{array}{c}
\lambda _1 \\
\rightarrow
\end{array}
\begin{array}{c}
\\
u_1
\end{array}
\begin{array}{c}
\lambda _2 \\
\rightarrow
\end{array}
\begin{array}{c}
\\
u_2
\end{array}
\begin{array}{c}
\\
\rightarrow
\end{array}
\begin{array}{c}
\\
...
\end{array}
\begin{array}{c}
\\
\rightarrow
\end{array}
\begin{array}{c}
\\
u_{s-1}
\end{array}
\begin{array}{c}
\lambda _s \\
\rightarrow
\end{array}
\begin{array}{c}
\\
v_q
\end{array}
,
\end{equation}
where $u_1,...,u_{s-1}$ are vertices of doubles (they may be absent, i.e. a
path may consist of exactly one arrow) and $z(\lambda _1)\not =0.$ Then by
Lemma 11 and item b) of the definition of a weight function, $z(\lambda
_1)=z(\lambda _2)=...=z(\lambda _s),$ this non-zero integer we shall call
the {\it weight} {\it of} {\it path} (2). We shall say that the elementary
path (2) {\it avoids} a vertex $v_i$ if $p<i<q.$ Now we establish some
properties of elementary paths:

A. The intersection of two elementary paths does not contain any vertex of a
double.

B. Each nondegenerate arrow is contained in an elementary path.

C. If a vertex $v_i$ is avoided by an elementary path (2) having length at
least 2, then the $v_i$ is incomparable with some vertex $u_l$ in this path.
Otherwise, $v_p<u_1<...<u_{s-1}<v_q$ implies one of the following
conditions: $v_p<v_i<u_1$ or $u_j<v_i<u_{j+1}$ for some $j$ or $%
u_{s-1}<v_i<v_q.$ This contradicts Lemma 9 because the vertices $%
u_1,...,u_{s-1}$ are contained in doubles.

D. If a vertex of a triple is avoided by an elementary path of length at
least 2, then all other vertices of this triple can not be avoided by any
elementary path having length $\geq 2.$ This follows from property C and
Lemma 8.

E. The sum of the weights of all elementary paths avoiding a vertex $v_i$ is
equal to $-z(v_i).$ Indeed, this is obvious for $v_1$ because, by property
B, only arrows having weight 0 can stop at $v_1.$ If property E is true for $%
v_i,$ then the sum of the weights of all elementary paths avoiding $v_i$ or
starting from $v_i$ is equal to 0. But the set of these paths coincides with
the set of all elementary paths avoiding $v_{i+1}$ or stopping at $v_{i+1}.$
Hence property E is true for $v_{i+1}.$

F. Let a triple $\{b_1,b_2,b_3\}$ satisfy the following conditions: 1) there
is no nondegenerate arrow starting from $a<b_1;$ 2) there is a pair of
connected degenerate strong arrows starting from $(b_1,b_2)$ or $(b_1,b_3);$
3) there is a pair of connected nondegenerate weak arrows starting from $%
(b_2,b_3).$ Then there exists a triple $\{a_1,a_2,a_3\}$ satisfying the same
conditions and $a_1<b_1.$ Indeed, let for definiteness the pair of connected
degenerate strong arrows start from $(b_1,b_2).$ From $z(b_1)=0,$ $%
z(b_2)=-z(b_3)\not =0$ and properties D and E, it follows that $b_2$ and $%
b_3 $ is avoided by a nondegenerate arrow. Let $b_3$ be avoided by a
nondegenerate arrow $a_i\rightarrow c_j.$ Then $a_i<b_3<c_j.$ By Lemma 9,
there exists a path $a_i\rightarrow ...\rightarrow b_3\rightarrow
...\rightarrow c_j$ consisting of strong arrows. But by Lemma 12, there is
only a weak arrow starting from $b_3.$ Hence $b_2$ is avoided by some
nondegenerate arrow $a_i\rightarrow c_j.$ By Lemma 9, it is a weak arrow, $%
i=j$ and there is a path $a_i\rightarrow ...\rightarrow b_2\rightarrow
...\rightarrow c_i$ consisting of strong arrows. But there is only one
strong arrow starting from $b_2$ and it is connected with an arrow starting
from $b_1.$ Hence the arrows connected with $a_i\rightarrow ...\rightarrow
b_2\rightarrow ...\rightarrow c_i$ compose the path $a_{i^{\prime
}}\rightarrow ...\rightarrow b_1\rightarrow ...\rightarrow c_{i^{\prime }}.$
The triple $\{a_1,a_2,a_3\}$ satisfies our requirement.

Let $c_l$ be the vertex such that there is a nondegenerate arrow starting
from $c_l$ and there is no nondegenerate arrow starting from $b<c_l$. Then
there is no nondegenerate arrow stopping at $c_l,$ hence $z(c_l)=0$ and
there are two arrows starting from $c_l$ and having the weights $n$ and $-n$%
, moreover, $l=1$ and the arrows connected with them start from $c_2$ and $%
c_3.$ Since $z(c_2)=-z(c_3)=\pm n\not =0,$ the vertices $c_2$ and $c_3$ are
avoided by elementary paths, and one of them is a nondegenerate arrow. Let
for definiteness $c_2$ be avoided by a nondegenerate arrow $b_i\rightarrow
d_j.$ By Lemma 9, $i=j$ and there is a path b$_i\rightarrow ...\rightarrow
c_2\rightarrow ...\rightarrow d_i.$ Since there exists exactly one arrow
starting from $c_2$ and this arrow is connected with an arrow starting from $%
c_1,$ we have that the arrows connected with b$_i\rightarrow ...\rightarrow
c_2\rightarrow ...\rightarrow d_i$ compose the path b$_{i^{\prime
}}\rightarrow ...\rightarrow c_1\rightarrow ...\rightarrow d_{i^{\prime }}.$
Since $b_{i^{\prime }}<c_1,$ there is no nondegenerate arrow starting from $%
b_{i^{\prime }}.$ Hence the arrow $b_i\rightarrow d_i$ is connected with the
arrow $b_{i^{\prime \prime }}\rightarrow d_{i^{\prime \prime }},$ where $%
i^{\prime }\not =i^{\prime \prime }$ and $i^{\prime }=1.$ By applying
property F to the triple $\{b_1,b_2,b_3\},$ we obtain a triple $%
\{a_1,a_2,a_3\}$. By applying property F to the triple $\{a_1,a_2,a_3\},$ we
obtain another triple and so on. This contradicts the finiteness of the
graph $\Gamma .$ This proves our Lemma.

{\bf Proof of Proposition 2.} We number all vertices and all arrows of the
graph$~\Gamma :$%
$$
\Gamma _0=\{a_1,a_2,...,a_r\},\qquad \Gamma
_1=\{f_{11},f_{12},...,f_{s1},f_{s2}\}.
$$
where $f_{j1}:a_{p(j1)}\rightarrow a_{q(j1)}$ and $f_{j2}:a_{p(j2)}%
\rightarrow a_{q(j2)}$ are two connected arrows and $a_{p(j1)}<a_{p(j2)}.$
Let the basis vector $m_i$ correspond to the vertex $a_i$ and let the double
morphism $f_i$ correspond to the pair $(f_{j1},f_{j2}).$ Then $%
f_jm_{p(j1)}=m_{q(j1)}$ and $f_jm_{p(j2)}=\lambda _jm_{q(j2)},$ where $%
\lambda _j$ is the parameter of a double morphism $f_j.$

By changes of the basis vectors
\begin{equation}
\label{d3}m_i=x_im_i^{\prime },\quad 0\not =x_i\in k,
\end{equation}
we obtain a new set of double morphism: $f_j^{\prime
}=x_{p(j1)}x_{q(j1)}^{-1}f_j,$ $1\leq j\leq s$, with the parameters $\lambda
_j^{\prime }=\lambda _jx_{p(j1)}x_{q(j1)}^{-1}x_{p(j2)}^{-1}x_{q(j2)}.$

The change (3) gives a multiplicative basis if $\lambda _1^{\prime }=\lambda
_2^{\prime }=...=\lambda _s^{\prime }=1,$ i.e. if $x_1,x_2,...,x_r$ satisfy
the system of equations
\begin{equation}
\label{d4}\lambda _jx_{p(j1)}x_{p(j2)}^{-1}=x_{q(j1)}x_{q(j2)}^{-1},\quad
1\leq j\leq s.
\end{equation}

We shall solve the system by elimination: solve the first equation for some $%
x_i$ and substitute the result in other equations. This amounts to the
multiplication of each of them by rational power of the first equation.
Further we solve the second equation of the obtained system for some $x_j$
and substitute the result in other equations... There are two possibilities:

1. After the $s$th step, we obtain the solution $(x_1,...,x_t)\in
(k\setminus \{0\})^t$ of (4).

2. After the $(t-1)$th step $(1\leq t\leq s)$, we obtain a system, the $t$th
equation of which does not contain unknowns. In this case, the $t$th
equation of (4), up to scalar multiples $\lambda _t$, is the product of
rational powers of the 1th,...,$(t-1)$th equations. It means that there
exist integers $z_1,...z_t$ such that $z_t\not =0$ and the equality
\begin{equation}
\label{d5}%
(x_{p(11)}x_{p(12)}^{-1})^{z_1}...(x_{p(t1)}x_{p(t2)}^{-1})^{z_t}=(x_{q(11)}x_{q(12)}^{-1})^{z_1}...(x_{q(t1)}x_{q(t2)}^{-1})^{z_t}
\end{equation}
is the identity, i.e. each $x_i$ has the same exponents at the two sides of
(5).

Define the integer function $z:\Gamma _1\rightarrow ${\cal Z} by $%
z(f_{j1})=-z(f_{j2})=z_j$ for $j\leq t$ and $z(f_{j1})=z(f_{j2})=0$ for $%
j>t. $ Since $x_i$ corresponds to the vertex $a_i$ of $\Gamma $, we have by
(5) that this function is a non-zero weight function, which contradicts
Lemma 13. Hence case 2 is impossible. This finishes the proof of Proposition
2.


\begin{thebibliography}{9}
\bibitem{BGRS}  Bautista R., Gabriel P., Roiter A.V., Salmeron L.
Representation-finite algebras and multiplicative basis. Invent. Math. {\bf %
81} (1985) 217--285.

\bibitem{GR}  Gabriel P., Roiter A.V. Representations of finite-dimensional
algebras. Encyclopaedia of Math. Sci., vol. 73, Algebra 8, Springer-Verlag
(1992), 177 p.

\bibitem{GNRSV}  Gabriel P., Nazarova L.A., Roiter A.V., Sergeichuk V.V.,
Vossieck D. Tame and wild subspace problems. Ukrainian Math. J. {\bf 45}
(no. 3) (1993) 313--352.
\end{thebibliography}
\end{document}